\documentclass[12pt]{amsproc}
\usepackage{amsmath,amssymb,euscript,Myapprox}
\newtheorem{theorem}{Theorem}[section]
\newtheorem{lemma}[theorem]{Lemma}
\newtheorem{definition}[theorem]{Definition}

\newtheorem{corollary}[theorem]{Corollary}

\newcommand{\eb}{\begin{equation}}
\newcommand{\ee}{\end{equation}}
\newcommand{\h}{Homeo(X)}

\newcommand{\e}{\varepsilon}
\newcommand{\aut}{Aut(X,\mathcal B)}

\begin{document}

\title[On approximation of homeomorphisms  of  a Cantor set]{On approximation of homeomorphisms  of \\  a Cantor set}

\date{}

\author{Konstantin Medynets}
\address{Department of Mathematics, Institute for Low Temperature Physics,
47 Lenin ave., 61103 Kharkov, Ukraine}
\email{medynets@ilt.kharkov.ua}

\maketitle
\begin{abstract} We continue to study topological properties of the group $Homeo(X)$ of
all homeomorphisms of a Cantor set $X$ with respect to the uniform
topology $\tau$, which was started in \cite{B-K 1}, [B-D-K 1; 2],
\cite{B-D-M}, and \cite{B-M}. We prove that the set of periodic
homeomorphisms is $\tau$-dense in $Homeo(X)$ and deduce from this
result that the topological group $(Homeo(X), \tau)$ has the
Rokhlin property, i.e., there exists a homeomorphism whose
conjugate class is $\tau$-dense in $Homeo(X)$. We also show that
for any homeomorphism $T$ the topological full group $[[T]]$ is
$\tau$-dense in the full group $[T]$.
\end{abstract}

%
%

\section{Introduction}
Many famous problems in ergodic theory involve the use of
topologies on the group $Aut(X,\mathcal B,\mu)$ of all
measure-preserving transformations of a standard measure space.
The first results on group topologies of $Aut(X,\mathcal B,\mu)$
appeared in the paper \cite{Hal 1}.  Halmos introduced two
topologies $d_u$ and $d_w$, which were called later the uniform
and weak topologies, respectively. He defined the uniform topology
$d_u$ by saying that two automorphisms $T$ and $S$ are ``close''
to each other if the quantity $\mu(\{x\in X : Tx\neq Sx\})$ is
small enough. The weak topology is generated by the  sets of the
form $N(T;E;\e)=\{S\in Aut(X,\mathcal B,\mu) : \mu(SE\triangle
TE)<\e\}$, where $T\in Aut(X,\mathcal B,\mu)$ and $E\in\mathcal
B$.  The use of these topologies turned out to be very fruitful
and led to many outstanding results in ergodic theory (for
references, see, for example, \cite{B-K-M} and \cite{C-F-S}). One
of the most relevant results in the theory is the Rokhlin lemma
\cite{Ro} stating that the set of periodic automorphisms
 is $d_u$-dense in $Aut(X,\mathcal
B,\mu)$.

The idea of investigation of transformation groups by means of
introducing various topologies into these groups was used in
\cite{B-D-K 1} and \cite{B-D-K 2}. In the papers, the authors
considered the groups $Aut(X,\mathcal B)$ of all automorphisms of
a standard Borel space and the group $\h$ of all homeomorphisms of
a Cantor set $X$ with respect to the several topologies analogous
to those in ergodic theory.

Following \cite{B-D-K 2}, we continue studying the group $\h$ of
all homeomorphisms of a Cantor set $X$ with respect to the
topology $\tau$ (cf. Definition \ref{top}), which is obviously a
direct analog of the topology $d_u$. We show that the set of all
periodic homeomorphisms is $\tau$-dense in $\h$ (Corollary
\ref{TopRokhlinLemma}). This result can be treated as a
topological version of the Rokhlin lemma. As a corollary, we prove
that the set of topologically free homeomorphisms is $\tau$-dense
in $\h$ (Theorem \ref{TopFreeInUniformTop}). Recall that a
homeomorphism is called {\it topologically free} if the set of
aperiodic points is dense.

In \cite{G-K},  an  interesting class of topological groups was
defined: by definition, a topological group has the {\it Rokhlin
property} if it has an element whose conjugate class is dense. The
authors raised the question: which groups possess this property.
At the moment, there is an extensive list of groups that have the
Rokhlin property. In particular, the group $\h$ with the topology
generated by the metric $D(T,S)=\sup_{x\in X}d(Tx,Sx)$, where $d$
is a metric on $X$ compatible with the topology \cite{Gl-W}, and
the group $Aut(X,\mathcal B,\mu)$ with respect to the weak
topology \cite{Hal 2} have the Rokhlin property. See also the
paper \cite{Ke-Ros} for a general approach to the study of groups
with dense conjugate classes.

Motivated by this, we present a unique approach allowing us to
show that the topological groups $(Aut(X,\mathcal B),\tau)$ and
$(\h,\tau)$ have the Rokhlin property (Theorem
\ref{RokhlinProperty}).

\medskip\noindent The other part of the paper is devoted to the study of
full groups $[T]$ of homeomorphisms $T\in \h$. Our motivation
comes from the paper  \cite{G-P-S}, where full groups were
indispensable in the study of orbit equivalence of Cantor minimal
systems. It is worthwhile to investigate full groups and their
dense subsets for any homeomorphism.

In this context, we show that  for any $T\in\h$, the topological
full group $[[T]]$ is $\tau$-dense in the full group $[T]$
(Theorem \ref{mainResult}).

In the last section, we give a description of homeomorphisms from
the topological full group $[[T]]$ of  aperiodic $T$ (Theorem
\ref{DescripFullGroups}). We consider a subgroup $\Gamma_Y$ of
$[[T]]$, which is an increasing union of permutation groups, and
find a criterion when $\Gamma_Y$  is $\tau$-dense in $[[T]]$
(Theorem \ref{gamma}).


\medskip \noindent {\bf Background.} Throughout the paper, $X$
denotes a {\it Cantor set} and $\mathcal B$ stands for the
$\sigma$-algebra of Borel sets of $X$. A one-to-one Borel map $T$
of $X$ onto itself is called an {\it automorphism} of $(X,\mathcal
B)$. Denote by $\aut$ the group of all automorphisms of
$(X,\mathcal B)$ and by $\h$ the group of all homeomorphisms of
$X$.

Following \cite{B-K 1}, recall  the definition of the {\it uniform
topology} on $\aut$. Let $\mathcal{M}_1(X)$ denote the set of all
Borel probability measures on $X$. For $T,S\in\aut$, denote
$E(T,S)=\{x\in X : Tx\neq Sx\}$.

\begin{definition}\label{top} The uniform topology $\tau$  on
$Aut(X,\mathcal{B})$ is defined by the base of neighborhoods
$\mathcal{U}= \{U(T;\mu_1,\ldots,\mu_n;\varepsilon)\}$,  where
$$\begin{array}{ll}
U(T;\mu_1,\mu_2,\ldots,\mu_n;\varepsilon)=\{S\in \aut\ |\
\mu_i(E(S,T))< \varepsilon,\;i=1,\ldots,n\}.
\end{array}
$$
Here $T\in\aut$, $\mu_1,\ldots,\mu_n\in\mathcal{M}_1(X)$, and
$\e>0$.
\end{definition}

As $\h$ is a subgroup of $\aut$, we also denote  by $\tau$ the
topology on $\h$ induced from $(Aut (X,\mathcal B),\tau)$.

Observe that $\aut$ and $\h$ are Hausdorff topological groups with
respect to the uniform topology $\tau$. More results related to
topological properties of $\aut$ and $\h$ with respect to $\tau$
can be found in \cite{B-D-K 1, B-D-K 2, B-D-M, B-K-M, B-M}.

Let $T\in\aut$. A point $x\in X$ is called {\it periodic} of
period $n>0$ if  $T^nx=x$ and $T^ix\neq x$ for $i=1,\ldots,n-1$.
If $T^nx\neq x$ for  $n\neq 0$, the point $x$ is called {\it
aperiodic}.  We say that  $T$ is  {\it aperiodic} if it has no
periodic points. Note that for any $T\in\aut$ the set $X$ can be
decomposed into a disjoint union of Borel sets
$X=X_\infty\cup\bigcup_{n\geq 1}X_n$, where $X_n$ consists of all
points of period $n$ and $X_\infty$ is formed by all aperiodic
points. Notice that some $X_n$'s can be empty. Moreover, for every
$X_n$, $n<\infty$, there exists a Borel set $X_n^0\subset X_n$
such that $X_n=\bigcup_{i=0}^{n-1}T^iX_n^0 $ is a disjoint union.
We call $\{X_\infty,X_1,X_2,\ldots\}$ the {\it canonical partition
of $X$ associated to} $T$.

 Recall that a finite
family of disjoint Borel sets $\xi=\{A,TA,\ldots,T^{n-1}A\}$ is
called a {\it $T$-tower} with the {\it base} $B(\xi)=A$  and the
{\it height} $h(\xi)=n$. A partition $\Xi=\{\xi_1,\xi_2,\ldots\}$
of $X$ is called a {\it Kakutani-Rokhlin} (K-R) partition if every
$\xi_i$ is a $T$-tower. For a K-R partition $\Xi$, we denote
$\bigcup_{n\geq 1}B(\xi_i)$ by $B(\Xi)$ and call it the {\it base
} of the K-R partition. Notice that for a K-R partition $\Xi$, one
has that $T^{-1}B(\Xi)=\bigcup_{\xi\in\Xi} T^{h(\xi)-1}B(\xi)$.

For $T\in\aut$, let $Orb_T(x)=\{T^nx: n\in\mathbb Z\}$ denote the
$T$-orbit of $x$. With any homeomorphism $T\in\h$, we can assign
two full groups $[T]_C$ and $[T]_B$, where
$$[T]_C=\{S\in\h : Orb_S(x)\subseteq
Orb_T(x),\,x\in X\}$$
 $$[T]_B=\{S\in\aut : Orb_S(x)\subseteq
Orb_T(x),\, x\in X\}.$$ Here the subindeces $C$ and $B$ stand for
the cases of Cantor and Borel dynamics, respectively. Clearly,
$[T]_C$ is a subgroup of $[T]_B$. Observe that if $S\in [T]_B$,
then there is a Borel function $n_S : X\rightarrow \mathbb Z$ such
that $Sx=T^{n_S(x)}x$ for all $x\in X$. The subgroup $[[T]]=\{S\in
[T]_C : n_S\mbox { is continuous}\}$  is called the {\it
topological full group} of $T$.

One of the main results in the approximation theory of Borel
automorphisms is a Borel version of the Rokhlin lemma. The
following $\tau$-version of the Rokhlin lemma was proved in
\cite[Proposition 3.6]{B-D-K 1}. We also refer the reader to the
works \cite[Section 7]{N} and \cite[Section 4]{W} for  measure
free versions of the result.

\begin{theorem}\label{RokhlinLemma} Let $T$ be an aperiodic  automorphism of
$X$. Then there exists a sequence of periodic automorphisms
$(P_n)\in Aut(X,\mathcal B)$  such that $P_n
\stackrel{\tau}{\longrightarrow} T,\ n\to \infty$. Moreover, the
automorphisms $P_n$ can be taken from $[T]_B$.
\end{theorem}

Denote by $\mathcal Per_0$ the set of all homeomorphisms $P$ such
that $P^n=\mathbb I$ for some $n\in\mathbb N$; and for $T\in\h$,
set $\mathcal Per_0(T)=\mathcal Per_0\cap [[T]]$.
%
%

\section{Rokhlin lemma}
In the section, we prove a topological version of the Rokhlin
lemma, namely, we show that the set of periodic homeomorphisms is
$\tau$-dense in $\h$. Then, we deduce several  corollaries of this
result. In particular, we prove that the topological group
$(\h,\tau)$ possesses the Rokhlin property and  the topological
full group $[[T]]$ is $\tau$-dense in $[T]_B$ for any $T\in\h$.

\begin{theorem}\label{mainResult} (1) The set $\mathcal Per_0$ is $\tau$-dense in $\h$.
\\
(2) Let $T\in\h$, then for any automorphism $S\in [T]_B$ and any
$\tau$-neighborhood $U=U(S;\mu_1,\ldots,\mu_p;\e)$ of $S$ there
exists a periodic homeomorphism $P\in [[T]]$ such that $P\in U$.
\end{theorem}
\textsf{Proof.} Notice that statement (1) is an immediate
corollary of (2).  By Theorem \ref{RokhlinLemma}, it is enough to
prove (2) for a periodic automorphism $S$.

Let us sketch the main stages of the proof. (i) We find a finite
number of disjoint $S$-towers consisting of closed sets and
covering  ``almost'' the entire space $X$ with respect to the
measures $\mu_i$ such that on each level of these $S$-towers the
automorphism $S$ coincides with a power of $T$. (ii) We extend the
$S$-towers found in (i) to clopen ones constructed by powers of
$T$. (iii) Using the clopen towers, we define a periodic
homeomorphism $P$ which belongs to $U(S;\mu_1,\ldots,\mu_p;\e)$.

(i) Let $\Xi=\{X_1,X_2,\ldots\}$ be the canonical Borel partition
of $X$ associated to $S$. Without loss of generality, we will
assume that the sets $X_i$ are not empty, $i\in\mathbb N$.

We first find $N\in\mathbb N$ such that
\begin{equation}\label{1}
\mu_j\left(X_1\cup X_2\cup\ldots\cup X_N
\right)>1-\frac{\e}{3}\mbox{ for all }j=1,\ldots,p.
\end{equation}
For $n\geq 1$, set $Z_n=\{x\in X : Sx=T^ix\mbox{ for some }-n\leq
i\leq n\}$. For $i\geq 1$ define the sets
$$X_i^0(n)=\bigcap_{j=0}^{i-1}S^{-j}\left(S^jX_i^0\cap Z_n
\right),$$ where $X_i=X_i^0\cup SX_i^0\cup\ldots\cup S^{i-1}X_i^0$
is a disjoint union.
 Since
$S\in[T]_B$, we have that $X_i=\bigcup_{n\geq 1}X_i(n)$, where
$X_i(n)=\bigcup_{j=0}^{i-1}S^jX_i^0(n)$. Then, find $K\in\mathbb
N$ such that
\begin{equation}\label{2}
\mu_j\left(\bigcup_{i=1}^N(X_i\setminus X_i(K))
\right)<\frac{\e}{3}\mbox{ for all }j=1,\ldots,p.
\end{equation}

Denote by $\mathcal S_i$  the set of all maps from
$\{0,\ldots,i-1\}$ to $\{-K,\ldots, K\}$. For $\sigma\in\mathcal
S_i$,  define the set
$$X_i^0(K,\sigma)=\bigcap_{j=0}^{i-1}S^{-j}\left(\{x\in S^jX_i^0(K) : Sx=T^{\sigma(j)}x \}\right).$$
Thus, we get a finite cover $X_i^0(K)=\bigcup_{\sigma\in\mathcal
S_i}X_i^0(K,\sigma)$. Applying the standard argument, make the
$X_i^0(K,\sigma)$'s disjoint and denote the obtained sets by
$X_i^0(K,\sigma)$ again. Some of the $X_i^0(K,\sigma)$'s can be
empty, but, without loss of generality, we will assume they are
not.  Observe that $S$ restricted to $S^jX_i^0(K,\sigma)$ is equal
 to $T^{\sigma(j)}$, for $i\geq 1$, $j=0,\ldots,i-1$, and
$\sigma\in\mathcal S_i$. This means that $S$ is a homeomorphism on
$S^jX_i^0(K,\sigma)$.

For every $X_i(K,\sigma)=\bigcup_{j=0}^{i-1}S^jX_i^0(K,\sigma)$,
find a closed set
 $A_i^0(\sigma)\subset X_i^0(K,\sigma)$
 such that
\begin{equation}\label{3}
 \mu_j(\bigcup_{i=1}^N\bigcup_{\sigma\in\mathcal S_i}(X_i(K,\sigma)\setminus A_i(\sigma)))<\frac{\e}{3}\mbox{ for }j=1,\ldots,p,
\end{equation}
where $A_i(\sigma)=\bigcup_{j=0}^{i-1}S^j A_i^0(\sigma)$.
%
%

(ii) Summing up the above, we get that $\{A_i(\sigma): 1\leq i
\leq N\;, \sigma\in\mathcal S_i\}$ is a family of disjoint closed
$S$-towers such that the automorphism $S$ restricted to
$S^jA_i^0(\sigma)$ is equal to $T^{\sigma(j)}$. Furthermore, it
follows from (\ref{1}),(\ref{2}), and (\ref{3}) that

\begin{equation}\label{4} \mu_j(\bigcup\limits_{i=1}^N\bigcup\limits_{\sigma\in\mathcal S_i}
A_i(\sigma))>1-\e.
\end{equation}

As the closed $S$-towers $A_i(\sigma)$'s  are disjoint, we can
 find clopen sets $\overline A_i^0(\sigma)\supset A_i^0(\sigma)$
 so that all the
 sets $\overline A_i^0(\sigma)$ and
  $T^{\sigma(0)+\ldots+\sigma(j)}\overline A_i^0(\sigma)$
  are mutually  disjoint  for $i=1,\ldots,N$ $j=0,\ldots,i-2$,
  and $\sigma\in \mathcal S_i$.

(iii) Define the periodic homeomorphism $P$ as follows:

$$Px=\left\{\begin{array}{ll}

\begin{array}{ll}
T^{\sigma(0)}x     &   \mbox{if }x\in \overline A_i^0 \\
T^{\sigma(j+1)}x & \mbox{if }x\in
T^{\sigma(0)+\ldots+\sigma(j)}\overline A_i^0(\sigma) \\
T^{-\sigma(0)-\ldots-\sigma(i-2)}x   &    \mbox{if }x\in
T^{\sigma(0)+\ldots+\sigma(i-2)}\overline A_i^0(\sigma)
\\
x & \textrm{otherwise}.
\end{array}  &

 \begin{array}{c} \sigma\in\mathcal S_i \\
     0\leq j\leq i-3
 \end{array}

\end{array}\right.$$

Clearly, $P$ is well-defined and belongs to $[[T]]$. By the
definition of $P$, we have
$$\{x\in X\ |\ Px=Sx\}\supset \bigcup_{i=1}^N\bigcup_{\sigma\in\mathcal S_i} A_i(\sigma).$$
Hence,  we get by (\ref{4}) that $P\in
U(S;\mu_1,\ldots,\mu_p;\e)$. This completes the proof.
\hfill$\square$ \medskip

\medskip\noindent {\bf Remark.} After this work was submitted, B.
Miller showed how using ideas of the proof above one can
generalize Statement (2) of Theorem \ref{mainResult} to any
countable group acting  by homeomorphisms on a zero-dimensional
Polish space \cite{Mil}.
%
%
%
\\
\\
{\bf Rokhlin property}
\\
\\
We give several immediate corollaries of Theorem \ref{mainResult},
which have the well-known analogs in ergodic theory.

\begin{corollary}\label{TopRokhlinLemma} Let $T\in\h$. Then, for
every $\tau$-neighborhood $U$ of $T$, there exists a homeomorphism
$P\in \mathcal Per_0(T)\cap U$ whose associated canonical
partition is clopen.
\end{corollary}

The next statement generalizes Theorem 4.5 of \cite{B-K 1} proved
originally for minimal homeomorphisms.

\begin{corollary}\label{denseInCanFullGroup} Let $T\in\h$. The topological full group $[[T]]$ of
$T$ is $\tau$-dense in $[T]_C$.
\end{corollary}

As the group $\h$ is not $\tau$-closed in $\aut$, in \cite{B-D-K
2} the authors brought up the question: how to describe the
closure of $[[T]]$ in $(\aut,\tau)$. They answered it for minimal
homeomorphisms (see Theorem 2.8 of \cite{B-D-K 2}) and we
generalize it up to an arbitrary homeomorphism.

\begin{corollary} Let $T\in\h$. Then
$\overline{[[T]]}^\tau=\overline{[T]_C}^\tau=[T]_B$.
\end{corollary}

%
%
\noindent{\bf Definition.} A topological group $G$ possesses the
{\it Rokhlin property} if the action of $G$ on itself by
conjugation is topologically transitive, i.e. there is an element
of $G$ whose conjugate class is dense.

The following proposition extends the list of topological groups
that have the Rokhlin property. See also \cite{Gl-W} and
\cite{Ke-Ros} for other examples.

\begin{theorem}\label{RokhlinProperty} The topological groups $(Aut(X,\mathcal B),\tau)$ and $(\h,\tau)$ possess the
Rokhlin property.
\end{theorem}
\textsf{Proof.} We prove this theorem for the group $(\h,\tau)$
only,  for  the  other case the proof is similar.

 Take a decomposition of the Cantor set
$X=\{x_0\}\cup\bigcup_{i\geq 1} X_i$ such that the $X_i$'s are
non-empty clopen sets with $diam(X_i\cup\{x_0\})\to 0$ as
$i\to\infty$. Let $S$ be a homeomorphism such that $Sx_0=x_0$ and
$S^ix=x$,  $S^jx\neq x$ for any $x\in X_i$, $j=1,\ldots,i-1$. Our
goal is to show that we can approximate any $T\in\h$ by elements
from the conjugate class of $S$. By Corollary
\ref{TopRokhlinLemma}, it suffices to approximate periodic
homeomorphisms whose canonical partitions are clopen. Thus,
suppose $T$ has a clopen partition $X=\bigcup_{i=1}^{k}Y_i$, where
$Y_i$ is the set of all points having $T$-period $n_i$ for some
$n_i\geq 0$. Observe that there exists a clopen set $Y_i^0$ such
that $Y_i=\bigcup_{j=0}^{n_i-1}T^jY_i^0$ is a disjoint union (see
Lemma 3.2 of \cite{B-D-K 2}). Analogously, there exists a clopen
set $X_{n_i}^0$ with $X_{n_i}=\bigcup_{j=0}^{n_i-1}S^jX_{n_i}^0$ a
disjoint union.

Let  $U=U(T;\mu_1,\ldots,\mu_p;\e)$ be a $\tau$-neighborhood of
$T$. Take a non-empty clopen $T$-invariant set $Z$ with
$\mu_i(Z)<\e$ for $i=1,\ldots,p$. Without loss of generality, we
may assume that $Y_i^0\setminus Z$ is not empty for
$i=1,\ldots,k$. Let $R_i$ be any homeomorphism from $X_{n_i}^0$
onto $Y_i^0\setminus Z$. Define a homeomorphism $R$ as follows:
let $R$ be equal to $T^jR_iS^{-j}$ whenever $x\in S^jX_{n_i}^0$
for $i=1,\ldots,k$, $j=0,\ldots,n_i-1$ and let $R$ map the rest of
the space $X$ onto $Z$. It is not hard to check that $RSR^{-1}\in
U$.

In the setting of Borel  dynamics, we need to produce a periodic
transformation that has uncountably many orbits of any finite
length. Then, the application of the Rokhlin lemma shows that its
conjugate class is dense. \hfill$\square$

\medskip\noindent {\bf Remark.} Let $p$ be the topology on $\h$ generated by the
metric $D(T,S)=\sup_{x\in X}d(Tx,Sx)$, where $d$ is a metric on
$X$ compatible with the topology. In \cite{Gl-W}, it is shown that
$(\h,p)$ has the Rokhlin property. Moreover, the elements whose
conjugate classes are dense form a residual set with respect to
$p$.
\\
\\
%
%
%
{\bf Topologically free homeomorphisms.}
\\
\\
It is interesting to compare the topological properties of the set
$\mathcal Ap$ of all aperiodic homeomorphisms with respect to the
both topologies $\tau$ and $p$. The following statement is proved
in \cite[Theorem 2.1]{B-D-K 2}.

\begin{theorem}\label{TopFreeInP} The set $\mathcal Ap$ is
dense in $(\h,p)$.
\end{theorem}

However, the situation in $(\h,\tau)$ is completely different. The
set $\mathcal Ap$ is nowhere dense with respect to the topology
$\tau $. To see this, one can check that the set $\mathcal Ap$ is
$\tau$-closed in $\h$. Then, the application of Theorem
\ref{mainResult} implies the result.

The question we investigate in this section is ``How can we extend
the class of aperiodic homeomorphisms to produce a $\tau$-dense
class?''. Apparently, the most natural extension of aperiodic
homeomorphisms is the class of topologically free homeomorphisms.
\\
\\
{\bf Definition.} It is said that a homeomorphism is {\it
topologically free} if the set of all aperiodic points is
dense.\medskip

\noindent  In Theorem \ref{TopFreeInUniformTop}, we prove that the
set of topologically free homeomorphisms is $\tau$-dense. To begin
with, we need the following lemma on homeomorphism extensions
proved in \cite{K-R}. We will need the arguments used in its
proof. Thus, we  give a sketch of the proof, but without going
into the details.

\begin{lemma}\label{extension} Let $A$ and $B$ be closed nowhere dense subsets of Cantor
sets $X$ and $Y$, respectively. Suppose there is a homeomorphism
$h: A\rightarrow B$. Then  $h$ can be extended to a homeomorphism
$h^* :X\rightarrow Y$ such that $h^*|_A=h$.
\end{lemma}
\textsf{Sketch of the proof.} Find clopen sets $\{U_i\}$ and
$\{V_j\}$ such that $X\setminus A=\bigsqcup_{i\geq 1} U_i$,
$Y\setminus B=\bigsqcup_{j\geq 1}V_j$, and their diameters tend to
zero. Find the points $a_i\in A$ such that
$dist(U_i,A)=dist(U_i,a_i)$ and $b_j\in B$ with
$dist(V_j,B)=dist(V_j,b_j)$.

Set $I=J=\mathbb N$. There exist injective functions $f:
I\rightarrow J$ and $g: J\rightarrow I$
  such that $$\begin{array}{ll}dist(U_i,a_i)>dist(V_{f(i)},h(a_i))\mbox{ for }i\in
  I\\
dist(V_j,b_j)>dist(U_{g(j)},h^{-1}(b_j))\mbox{ for }j\in J.
  \end{array}$$
Applying the usual Schr\"{o}der-Bernshtein argument to $f$ and
$g$, find disjoint sets $I=I'\sqcup I''$ and $J=J'\sqcup J''$ such
that $f(I')=J'$ and $g(J'')=I''$.

Let $\phi$ be an arbitrary homeomorphism of $U'=\bigcup_{i\in
I'}U_i$ onto $V'=\bigcup_{j\in J'}V_i$ such that
$\phi(U_i)=V_{f(i)}$. Analogously, let $\psi$ be a homeomorphism
of $V''=\bigcup_{j\in J''}V_j$ onto $U''=\bigcup_{i\in I''}U_i$
such that $\psi(V_j)=U_{g(j)}$.

 Define $$h^*(x)=\left\{
\begin{array}{ccc}  \phi(x) & x\in
U'\\
\psi^{-1}(x ) & x\in U''\\
h(x) & x\in A.\end{array}\right.$$ For the verification of
continuity of $h^*$, we refer the reader to \cite{K-R}.
\hfill$\square$

\begin{theorem}\label{TopFreeInUniformTop}  The set of topologically
free homeomorphisms is $\tau$-dense in $\h$.
\end{theorem}
\textsf{Proof.} By Corollary \ref{TopRokhlinLemma}, it suffices to
approximate only homeomorphisms from $\mathcal Per_0$ that have
clopen canonical partitions. Assume that $R$ belongs to $\mathcal
Per_0$ and its canonical partition $X=X_{n_1}\cup\ldots\cup
X_{n_m}$ is clopen. Recall that the set $X_{n_i}$ consists of all
points with the period $n_i$. Consider a $\tau$-neighborhood
$U=U(R;\mu_1,\ldots,\mu_k;\e)$. Since the $X_{n_i}$'s are
$R$-invariant and clopen, we will prove the theorem under the
assumption that  $X=X_{n_i}$ for some $i$ and leave the
generalization to the reader.

 Suppose $X=\bigcup_{i=0}^{p-1}R^iF$
is a clopen partition and $R^px=x$ for all $x\in X$. Using the
standard Cantor argument, find a closed nowhere dense set
$P\subset R^{p-1}F$ such that $\mu_i(P)>1-\e$ for $i=1,\ldots,k$.
Repeating the proof of Lemma \ref{extension}, we extend the
homeomorphism $R:P\rightarrow RP$ to a homeomorphism $T:
R^{p-1}F\rightarrow F$ so that the homeomorphism $T^*\in\h$
defined as $T^*|_{R^{p-1}F}=T|_{R^{p-1}F}$ and $T^*=P$ elsewhere
is topologically free. To do this, it suffices to choose the
functions $\psi$ and $\phi$ so that $\phi(x)\neq Rx$ and
$\psi^{-1}(x)\neq Rx$ for $x\in R^{p-1}F\setminus P$. Since
$E(T^*,R)=R^{p-1}F\setminus P$, we get that $T^*\in U$.
\hfill$\square$

\section{Structure of homeomorphisms from topological full group}
In this section, we discus the  structure of homeomorphisms from
the topological full group $[[T]]$ for  arbitrary aperiodic
$T\in\h$.

Consider a Cantor aperiodic system $(X,T)$. A Borel set $Y\subset
X$ is called {\it wandering} if $T^nY\cap Y=\emptyset$ for all
$n\geq 1$.

 {\bf Definition.} W say that a closed wandering set $Y$ is
 {\it basic} if every clopen neighborhood of $Y$ meets
 every  $T$-orbit.

\begin{theorem} Every Cantor aperiodic system has a basic set.
\end{theorem}
\textsf{Sketch of the proof.} Applying the argument developed in
\cite[Theorem 2]{B-D-M}, we can find a decreasing sequence of
clopen sets $\{U_n\}$ such that:  $U_{n+1}\subset U_n$;
$T^iU_n\cap U_n=\emptyset$ for $i=1,\ldots,n-1$; and $U_n$ meets
every $T$-orbit. Then $Y=\bigcap_n U_n$ is a  basic set.
\hfill$\square$\medskip

\noindent{\bf Remark} For  more results related to basic sets and
their interaction with Bratteli diagrams, see \cite{M}.\medskip

\noindent  Fix a triple $(X,T,Y)$, where $(X,T)$ is a Cantor
aperiodic system and $Y$ is a basic set. Consider a clopen
neighborhood $U$ of $Y$. It is not hard to check that for every
$x\in U$, there is $n=n(x)>0$ such that $T^nx\in U$. Therefore, it
follows from the definition of a basic set that, by the first
return function, we can construct a clopen K-R partition $\Xi$ of
$X$ with the base $B(\Xi)=U$.

Take a decreasing sequence of clopen sets $\{U_n\}$ such that
$Y=\bigcap_n U_n$. Constructing clopen K-R partitions for the
$U_n$'s and refining them, we  prove the following:

\begin{theorem}\label{KR-Partitions} Let $(X,T,Y)$ be a Cantor aperiodic system with a basic set $Y$. There
exists a sequence of clopen K-R partitions $\{\mathcal P_n\}$  of
$X$ such that for all $n\geq 1$
 (i) $\mathcal P_{n+1}$ refines $\mathcal P_n$;  (ii) $h_{n+1}>h_n$,
  where $h_n$ is the height of the lowest $T$-tower in
 $\mathcal P_n $; (iii) $B(\mathcal P_n)\supset B(\mathcal
P_{n+1})$; (iv) the sequence $\{\mathcal P_n\}$ generates the
clopen topology of $X$; (v) $\bigcap_nB(\mathcal P_n)=Y$.
\end{theorem}

We will follow here the method developed in \cite{B-K 1} for
minimal homeomorphisms (see also \cite{K-W}). Let $\mathcal P$ be
a clopen K-R partition with towers $\mathcal P(i)$,
$i=1,\ldots,k$. Define two partitions $\alpha=\alpha(\mathcal P)$
and $\alpha'=\alpha'(\mathcal P)$ of $\{1,2,\ldots,k\}$. We say
that $J$ is an atom of $\alpha $ if there exists a subset $J'$ of
$\alpha'$ such that
\begin{equation}\label{atom}
T(\bigcup_{i\in J}T^{h(i)-1}D_i)=\bigcup_{i'\in J'}D_{i'}
\end{equation}
and for every  proper subset $J_0$ of $J$, the $T$-image of
$\bigcup_{i\in J_0}T^{h(i)-1}D_i$ is not  a union of atoms from
$\mathcal P$. It follows from (\ref{atom}) that $J'$ is uniquely
defined by $J$ and $T$.

Let $S\in [[T]]$. Then,  there are a finite set $K\subset \mathbb
Z$ and clopen partition $\mathcal{E}=\{E_k : k\in K\}$ of $X$ such
that $Sx=T^kx$ for $x\in E_k$ and $k\in K$. Denote by $\mathcal
E(K)$ the clopen partition $\{S^kE_k : k\in K\}$. By Theorem
\ref{KR-Partitions}, find a K-R partition $\mathcal P=\{\mathcal
P(i) : i=1,\ldots,k\}$ with $\mathcal
P(i)=\{D_{0,i},\ldots,D_{h(i)-1,i}\}$ and $D_{j+1,i}=TD_{j,i}$
that refines $\mathcal E$ and $\mathcal E(K)$ and so that
$K\subset (-h,h)$, where $h$ is the height of the lowest $T$-tower
in $\mathcal P$.

Let $\mathcal F=\{(j,i) | i=1,\ldots,k,\,j=0,\ldots,h(i)-1\}$.
Observe that for every pair $(j,i)\in \mathcal F$ there is a
unique $l=l(j,i)\in K$ such that
\begin{equation}\label{power}S(D_{j,i})=T^lD_{j,i}.\end{equation}

Divide $\mathcal F=\mathcal F(\mathcal P)$  into three disjoint
sets $\mathcal F_{in}, \mathcal F_{top}$ and $\mathcal F_{bot}$ as
follows:
\begin{enumerate}
\item[(a)] $(j,i)\in \mathcal F_{in}$ if $S(D_{j,i})\subset
\mathcal P(i)$, i.e. $0\leq l+j\leq h(i)-1$; \item[(b)] $(j,i)\in
\mathcal F_{top}$ if $S(D_{j,i})$ goes through the top of
$\mathcal P(i)$, i.e. $l+j\geq h(i)$; \item[(c)] $(j,i)\in
\mathcal F_{bot}$ if $S(D_{j,i})$ goes through the bottom of
$\mathcal P(i)$, i.e. $l+j<0$,
\end{enumerate}
here $l$ is taken from (\ref{power}).

Let $\alpha$ and $\alpha'$ be the partitions of $\{1,\ldots,k\}$
defined by $T$ and $\mathcal P$. For  $J\subset \{1,\ldots, k\}$,
set $h_J=\min\{h(i) | i\in J\}$. For $J\in \alpha$ and
$J'\in\alpha'$, let $$F_1(r,J)=\bigcup_{i\in
J}D_{h(i)-h_J+r,i}\qquad F_2(r',J')=\bigcup_{i'\in J'}D_{r,i'}$$
where $r=0,\ldots,h_{J}-1$ and $r'=0,\ldots,h_{J'}-1$.\medskip

\noindent{\bf Definition.} We say that $S\in [[T]]$ belongs to
$\Gamma(\mathcal P)$ if for each pair $(j,i)\in \mathcal F$ the
following conditions hold:

\begin{enumerate} \item[(a)] if $(j,i)\in \mathcal F_{top}$ and $D_{j,i}\subset
E_l$, then $F_1(h_J-h(i)+j,J)\subset E_l$, where $J$ is an atom of
$\alpha$ containing $i$;

\item[(b)] if $(j,i)\in \mathcal F_{bot}$ and $D_{j,i}\subset
E_l$, then $F_2(j,J')\subset E_l$, where $J'$ is an atom of
$\alpha'$ containing $i$.
\end{enumerate}

 Condition (a) means that whenever the set $D_{j,i}$ goes through
the top of $\mathcal P(i)$ under the action of $S$, then the
entire level $F_1(r,J)$  containing $D_{j,i}$ also goes through
the top of $\mathcal P$. Similarly, one can clarify condition (b)
by taking  the $D_{j,i}$'s and levels $F_2(j,J')$ containing them
that go through the bottom of $\mathcal P$. Observe that if
$(j,i)\in \mathcal F_{in}$, then the entire levels $F_1(r,J)$ and
$F_2(j,J')$ containing $D_{j,i}$ remain ``within'' $\mathcal
P$.\medskip

Clearly, $\Gamma(\mathcal P)$ is a finite set. The following
theorem reveals the structure of homeomorphisms from $[[T]]$ for
an arbitrary aperiodic homeomorphism $T$. Notice that this
structure was found earlier for minimal homeomorphisms (see
Theorem 2.2 in \cite{B-K 1}).  Since  our proof is similar to that
in \cite[Theorem 2.2]{B-K 1}, we omit it.

\begin{theorem}\label{DescripFullGroups} Let $(X,T,Y)$ be a Cantor aperiodic system with a basic set $Y$ and  a
sequence of K-R partitions $\{\mathcal P_n\}$ satisfy the
conditions of Theorem \ref{KR-Partitions}. Then
$[[T]]=\bigcup_n\Gamma(\mathcal P_n)$ with $\Gamma(\mathcal
P_n)\subset \Gamma(\mathcal P_{n+1})$.
\end{theorem}
%

\noindent{\bf The subgroup $\Gamma_Y$}

Let $(X,T)$  be a Cantor aperiodic system with a basic set $Y$.
Define the subgroup $\Gamma_Y$ of $[[T]]$ as follows:
$S\in\Gamma_Y$ if $S\in \Gamma(\mathcal P_n)$ (and hence $S\in
 \Gamma(\mathcal P_m)$, for $m>n$) implies that $\mathcal F(\mathcal P_n)=\mathcal F_{in}$. In other
words, $S\in \Gamma_Y$ if no level from $\mathcal P_n$ goes over
the top as well as through the bottom under the action of $S$.
This means that $S$ acts as a permutation on each $T$-tower from
$\mathcal P_n$. Therefore, the group $\Gamma_Y$ is an increasing
union of permutation groups.

Our object is to find a criterion when $\Gamma_Y$ is dense in
$[T]$.

\noindent\medskip{\bf Remark.} (1) Denote by $[[T]]_Y$ the
subgroup of [[T]] consisting of homeomorphisms that preserve the
forward $T$-orbit of every $y\in Y$, i.e., $S\in [[T]]_Y$ if
$S(\{T^ny : n\geq 0\})=\{T^ny : n\geq 0\}$.  Observe that
$\Gamma_Y\subset [[T]]_Y$.

(2) The subgroup $[[T]]_Y$ is not $\tau$-dense in $[T]$.
Therefore, so is $\Gamma_Y$. Indeed, take any $z\in T^{-1}Y$ and
the Dirac measure $\delta_z$ supported by $\{z\}$. Consider  $S\in
U:=U(T;\delta_z;1/2)$. As $z\notin Y$ and $Sz=Tz\in Y$,  $S$ does
not preserve the forward $T$-orbit of $Tz$. Therefore, $U$
contains no elements from $[[T]]_Y$.
\\
\\
The fact that $\Gamma_Y$ is not $\tau$-dense in $[[T]]$ is mainly
caused by the presence of discrete measures. We can partly
overcome this obstacle by considering only continuous measures  in
the definition of the topology $\tau$. Denote by $\tau_0$ the
topology defined by continuous measures as in Definition
\ref{top}. One can check that $\tau_0$ is a Hausdorff group
topology on $\h$. The next theorem answers the question when
$\Gamma_Y$ is $\tau_0$-dense in $[T]$.

\begin{theorem}\label{gamma} Suppose we have a Cantor aperiodic system $(X,T)$
with a basic set $Y$. Then the subgroup $\Gamma_Y$ is
$\tau_0$-dense in $[T]$ if and only if the basic set $Y$ is at
most countable.
\end{theorem}
\textsf{Proof.} (1) Assume that $Y$ is uncountable. Take any
continuous measure $\mu$ supported by $T^{-1}Y$. Then for every
$S\in U:=U(T;\mu;1/2)$ there is at least one $z\in T^{-1}Y$ such
that $Sz=Tz$. This implies that $\{T^n(Tz) : n\geq 0\}$ is not
$S$-invariant. Therefore, by (1) of the remark above we get that
$\Gamma_Y\cap U=\emptyset$.

(2) Now, assume that  $Y$ is countable. Observe that by Corollary
\ref{denseInCanFullGroup}, it is enough to approximate
homeomorphisms from $[[T]]$ with elements of $\Gamma_Y$. Consider
$R\in [[T]]$ and a $\tau_0$-neighborhood
$U=U(R;\mu_1,\ldots,\mu_p;\e)$ of $R$. By  definition of $R$, the
sets $E_k=\{x\in X : Rx=T^kx\}$, $k\in K$, $|K|<\infty$, form a
clopen partition of $X$. Let $k_0=\sup \{|k| : k\in K\}$. As $Y$
is countable,  $\mu(T^nY)=0$ for any continuous measure $\mu$ and
integer $n$. Therefore, by Theorem \ref{KR-Partitions} we can find
a K-R partition $\mathcal P_n$ such that $k_0<2h_n$, where $h_n$
is the height of the lowest $T$-tower in $\mathcal P_n$, and

\begin{equation}\label{last}\mu_j(\bigcup\limits_{i=-k_0}^{k_0} T^iB(\mathcal
P_n))<\e,\mbox{ for }j=1,\ldots,p.\end{equation}
 Define a
homeomorphism $S\in\Gamma_Y\cap U$  as follows:

Take a $T$-tower (say $\lambda=\{D,\ldots,T^h(\lambda) D\}$) from
$\mathcal P_n$. Consider an atom  $T^lD$ of $\lambda$. We have two
possibilities:

(i) The $R$-orbit of the set $T^{l}D$ does not leave the $T$-tower
$\lambda$. In this case, we define the homeomorphism $S$ to be
equal to $R$ on the $R$-orbit of $T^{l}D$.

(ii) The set $T^{l}D$ leaves $\lambda$ under the action of $R$.
Then there exist integers $q<0<d$ such that  $R^{d+1}T^{l}D$ and
$R^{q-1}T^{l}D$ do not lie in $\lambda$ entirely, whereas the sets
$R^jT^{l}D$, $j=q,\ldots,d$ are contained in $\lambda$.   In this
case, we set $S=R$ on $R^iT^{l}D$, $i=q,\ldots,d-1$, and
$S=R^{-d+q}$ on $R^dT^{l}D$.

Observe that the choice of $\mathcal P_n$ guarantees  that
$$R^dT^l D\subset \bigcup\limits_{i=0}^{k_0} T^iD\cup\bigcup\limits_{i=h(\lambda)-1-k_0}^{h(\lambda)-1} T^iD.$$

Clearly, the homeomorphism $S$ constructed in (i) and (ii) is
periodic, but it is not defined yet on the entire space. To expand
its domain, we consider an atom $T^wD$ of $\lambda$ on which $S$
is not defined yet, and repeat (i) and (ii) with $T^wD$.

Repeating this procedure with every atom of $\lambda$, we define
$S$ on  $\lambda$. Moreover, $\lambda$ is $S$-invariant. By
construction, $S$ coincides with $P$ everywhere, maybe, except for
the set $\bigcup_{i=0}^{k_0}
T^iD\cup\bigcup_{i=h(\lambda)-1-k_0}^{h(\lambda)-1} T^iD.$

To finish constructing $S$, we need to repeat the argument for
every $T$-tower of $\mathcal P_n$.

The definition of $S$ implies that $S\in\Gamma_Y$ and $$\{x\in X :
Sx\neq Rx\}\subset \bigcup\limits_{i=-k_0}^{k_0} T^iB(\mathcal
P_n).$$

Therefore, by (\ref{last}) we have that $S\in U$. \hfill$\square$
\\
\\
{\it Acknowledgement.} This work was done when the author was
visiting the Torun University. I am thankful to the institution
for the hospitality and support. Also I would like to thank S.
Bezuglyi and J. Kwiatkowski for helpful discussions.

%


\begin{thebibliography}{zzz}

\bibitem[B-D-K 1]{B-D-K 1} S.~Bezuglyi, A.H.~Dooley, and J.~Kwiatkowski,
Topologies on the group of Borel automorphisms of a standard Borel
space, {\it Topol. Methods in Nonlinear Anal.}, {\bf 27} (2006),
333-385.

\bibitem[B-D-K 2]{B-D-K 2} S.~Bezuglyi, A.H.~Dooley, and J.~Kwiatkowski,
Topologies on the group of homeomorphisms of a Cantor set, {\it
Topol. Methods in Nonlinear Anal.}, {\bf 27} (2006), 299-331.


\bibitem[B-D-M]{B-D-M} S.~Bezuglyi, A.H.~Dooley, and K.~Medynets, The
Rokhlin lemma for homeomorphisms of a Cantor set, {\it Proc. Amer.
Math. Soc.} 133 (2005), 2957-2964.

\bibitem[B-K]{B-K 1} S.~Bezuglyi, J.~Kwiatkowski, The topological full
group of a Cantor minimal system is dense in the full group, {\it
Topol. Methods in Nonlinear Anal.}, {\bf 16} (2000), 371 - 397.


\bibitem[B-K-M]{B-K-M} S.~Bezuglyi, J.~Kwiatkowski, and
K.~Medynets,   Approximaion in measurable, Borel, and Cantor
dynamics, {\it Contemp. Math. Amer. Math. Society}, Volume 385,
2005.

\bibitem[B-M]{B-M} S.~Bezuglyi, K.~Medynets, Smooth automorphisms
and path-connectedness in Borel dynamics, {\it Indag. Math.} {\bf
15}, no. 4, (2004), 453-468.

\bibitem[C-F-S]{C-F-S} I.~Cornfeld, S.~Fomin, Ya.~Sinai, {\it Ergodic
Theory}, Grundlehren der mathematischen Wissenschaften {\bf 245},
Springer-Verlag, 1982.


\bibitem[G-K]{G-K} E. Glasner and J. King, A zero-one law for
dynamical properties,  Topological Dynamics and Applications (A
volume in honor of Robert Ellis), {\it Contemp. Math.}, vol. 215,
AMS, 1998, 215-242.


\bibitem[G-P-S]{G-P-S} T.~Giordano, I.~Putnam, and C.~Skau, Full
groups of Cantor minimal systems, {\it Israel. J. Math.}, {\bf
111} (1999), 285 - 320.


\bibitem[Gl-W]{Gl-W} E.~Glasner and B.~Weiss,  The topological Rohlin
property and topological entropy, {\it Amer. J.  Math.}, {\bf 123}
(2001), 1055 - 1070.

\bibitem[Hal 1]{Hal 1} P. Halmos, Approximation theories for
measure-preserving transformations, {\it Trans. Amer. Math. Soc.},
{\bf 55} (1944), 1-18.

\bibitem[Hal 2]{Hal 2} P.~Halmos,  {\it Lectures on Ergodic Theory},
 The Mathematical Society of Japan, Publications of the Mathematical
 Society of Japan, no. 3, 1956, vii+99 pp.

\bibitem[Kn-R]{K-R} B.~Knaster and M.~Reichbach, Notion
d'homog\'en\'eit\'e des hom\'eomorphies, {\it Fund. Math.}, {\bf
40} (1953), 180-193.


\bibitem[Ke-Ros]{Ke-Ros} A. Kechris, C. Rosendal, Turbulence, amalgamation
and generic automorphisms of homogenuous structures,
ArXiv:math.LO/0409567 v2, 30 September 2004.

\bibitem[K-W]{K-W} J. Kwiatkowski,  M. Wata,
Dimension and infinitesimal groups of Cantor minimal systems {\it
Topol. Methods Nonlinear Anal.}, 23, No.1,  (2004), 161-202.

\bibitem[M]{M} K. Medynets, Cantor aperiodic systems and Bratteli
diagrams, {\it C. R., Math., Acad. Sci. Paris}, {\bf 342}, No. 1,
(2006), 43-46.

\bibitem[Mil]{Mil} B. Miller, Density of topological full groups, preprint, 2006.


\bibitem[N]{N} M.~Nadkarni, {\it Basic Ergodic Theory}, 2nd Edition,
Birkh\"auser, 1998.

\bibitem[Ro]{Ro} V.A.~Rokhlin,  Selected topics from the metric theory
of dynamical systems (Russian),  {\it Usp. Mat. Nauk (N.S.)}, {\bf
4} (1949), no. 2, 57--128 (Engl. Transl. in `Amer. Math. Soc.
Translations',
 {\bf 49} (1966), 171 -- 240).

\bibitem[W]{W} B. Weiss, Measurable dynamics, {\it Contemp.
Math.}, {\bf 26} (1984), 395-421.


\end{thebibliography}
\end{document}